\documentclass{amsart}

\def\eps {{\epsilon}}

\begin{document}

\title[Triangle Inequality for Wasserstein Distances]{A Duality-Based Proof of the Triangle Inequality\\ for the Wasserstein Distances}

\author{Fran\c cois Golse}
\address{\'Ecole polytechnique \& IP Paris, CMLS, 91128 Palaiseau Cedex, France}
\email{francois.golse@polytechnique.edu}

\begin{abstract}
This short note gives a proof of the triangle inequality based on the Kantorovich duality formula for the Wasserstein distances of exponent $p\in[1,+\infty)$ in the case of a general Polish space. 
In particular it avoids the ``glueing of couplings'' procedure used in most textbooks on optimal transport. 
\end{abstract}



\maketitle




\section{Introduction}


Let $(\mathcal E,d)$ be a Polish metric space (with metric denoted by $d$). The Wasserstein distances are a family of metrics defined on subsets of $\mathcal P(\mathcal E)$, the set of Borel probability 
measures on $\mathcal E$. Specifically, for each $p\ge 1$, set
\[
\mathcal P_p(\mathcal E):=\left\{\mu\in\mathcal P(\mathcal E)\text{ s.t. }\int_\mathcal Ed(x_0,x)^p\mu(dx)<\infty\right\}\,.
\]
(Obviously, if the above condition is satisfied for one $x_0\in\mathcal E$, it is satisfied for all $x_0\in\mathcal E$ by the triangle inequality for the metric $d$ and by convexity of the map $z\mapsto z^p$
on $(0,+\infty)$.) For all $\mu,\nu\in\mathcal P(\mathcal E)$, a coupling of $\mu,\nu$ is an element $\rho$ of $\mathcal P(\mathcal E\times\mathcal E)$ such that\footnote{Since $\mathcal E$ is separable,
the Borel $\sigma$-algebra of $\mathcal E\times\mathcal E$ is the product of the Borel $\sigma$-algebra of $\mathcal E$ with itself: see Proposition 2.4.2 in chapter I of \cite{Mallia}.}
\[
\iint_{\mathcal E\times\mathcal E}(\phi(x)+\psi(y))\rho(dxdy)=\int_\mathcal E\phi(x)\mu(dx)+\int_\mathcal E\psi(y)\nu(dy)
\]
for all $\phi,\psi\in C_b(\mathcal E)$ (the set of bounded continuous functions on $\mathcal E$). The set of couplings of $\mu,\nu$ is denoted by $\mathcal C(\mu,\nu)$. (Notice that $\mathcal C(\mu,\nu)$
is nonempty since $\mu\otimes\nu\in\mathcal C(\mu,\nu)$.) The Wasserstein distance with exponent $p$ between $\mu,\nu\in\mathcal P_p(\mathcal E)$ is defined as
\[
\mathcal W_p(\mu,\nu):=\left(\inf_{\rho\in\mathcal C(\mu,\nu)}\iint_{\mathcal E\times\mathcal E}d(x,y)^p\rho(dxdy)\right)^{1/p}<+\infty\,.
\]
The Wasserstein distances are of considerable importance in the field of optimal transport, with applications to the calculus of variations, to statistics, to statistical mechanics, to machine learning, to cite
only a few applications --- see \cite{Rachev2} for more examples. For instance, $\mathcal W_p$ metrizes a variant of the weak convergence of probability measures in $\mathcal P_p(\mathcal E)$: see 
Theorem 7.12 in \cite{VillaniAMS}, or Theorem 6.9 in \cite{VillaniTOT}.

Recently, analogues of the Wasserstein metric with exponent $2$ have been defined in the quantum setting \cite{FGMouPaul,PalmaTrev}. Specifically, these analogues measure the difference between
two density operators, i.e. nonnegative self-adjoint trace-class operators with trace equal to one defined on a separable Hilbert space, which are the quantum analogue of phase-space probability measures
in classical mechanics. These analogues of the Wasserstein distance of exponent $2$ satisfy some variant of the triangle inequality (see Theorem 5.1 in \cite{FGPaulJFA} and formula (51) in \cite{PalmaTrev}) 
but, at the time of this writing, whether the genuine triangle inequality is satisfied by these quantum Wasserstein (pseudo)metrics remains an open question. 

If one returns to the classical setting, most references on optimal transport prove the triangle inequality for the Wasserstein distances by means of a procedure known as ``glueing couplings'' between
Borel probability measures. Specifically, given $\lambda,\mu,\nu\in\mathcal P_p(\mathcal E)$, pick $\rho_{12}\in\mathcal C(\lambda,\mu)$ and $\rho_{23}\in\mathcal C(\mu,\nu)$; the glueing procedure
provides us with $\sigma\in\mathcal P(\mathcal E\times\mathcal E\times\mathcal E)$ such that
\[
\begin{aligned}
\iiint_{\mathcal E\times\mathcal E\times\mathcal E}(\Phi(x,y)+\Psi(y,z))\sigma(dxdydz)
\\
=\iint_{\mathcal E\times\mathcal E}\Phi(x,y)\rho_{12}(dxdy)+\iint_{\mathcal E\times\mathcal E}\Psi(y,z)\rho_{23}(dydz)
\end{aligned}
\]
for all $\Phi,\Psi\in C_b(\mathcal E\times\mathcal E)$. This is the key step in one proof of the triangle inequality for $\mathcal W_p$ (see for instance \cite{VillaniAMS,AmbroGigliSava,VillaniTOT,Figalli}).
Once the probability measure $\sigma$ has been constructed, the remaining part of the proof of the triangle inequality is a routine computation involving the Minkowski inequality
One proof of the existence of $\sigma$ provided by the glueing procedure is based on disintegration of probability measures (see for instance Lemma 7.6 in \cite{VillaniAMS} or Remark 5.3.3 in 
\cite{AmbroGigliSava}). There exists another argument avoiding disintegration of measures, which is based on the Hahn-Banach theorem in the special case where $\mathcal E$ is compact: see 
Exercise 7.9 in \cite{VillaniAMS}. Still another proof of the triangle inequality uses an optimal transport map, when it exists: see \cite{Santambro}. An optimal transport map between the probability
measures $\mu,\nu\in\mathcal P_p(\mathcal E)$ is a Borel measurable map $T:\,\mathcal E\to\mathcal E$ such that\footnote{We denote by $T_\#\mu$ the image of the probability measure $\mu$
by the map $T$, defined by 
\[
\int_\mathcal E\phi(y)T_\#\mu(dy)=\int_\mathcal E\phi(T(x))\mu(dx)\,,\quad\text{ for all }\phi\in C_b(\mathcal E)\,.
\]}
\[
T_\#\mu=\nu\quad\text{ and }\quad\mathcal W_p(\mu,\nu)=\left(\int_\mathcal E d(x,T(x))^p\mu(dx)\right)^{1/p}
\]
The existence and uniqueness of an optimal transport map that is the gradient of a convex function is known in the case where $\mathcal E$ is a (finite-dimensional) Euclidean space, with Euclidean
metric $d$, and $\mu$ is absolutely continuous\footnote{In fact, it is enough to assume that $\mu(A)=0$ for all Borel measurable $A\subset\mathcal E$ of Hausdorff codimension $\ge 1$.} with respect 
to the Lebesgue measure on $\mathcal E$: this is Brenier's theorem (see Theorem 2.12 (ii) in \cite{VillaniAMS}).

None of the ingredients mentioned above (existence of optimal transport maps, glueing of couplings) are known to have analogues in general in the quantum setting. In fact, a recent counterexample
due to D. Serre \cite{DS} shows that the glueing procedure \textit{cannot} be extended to the quantum setting for \textit{arbitrary} couplings. 

In view of all these obstructions, we propose still another proof of the triangle inequality for the Wasserstein distances, in the hope that a different approach could perhaps lead to a better understanding 
of the quantum case.


\section{Kantorovich Duality}


First, we recall the Kantorovich duality for the Wasserstein distance $\mathcal W_p$.

\noindent
\textbf{Kantorovich Duality Theorem.} Let $(\mathcal E, d)$ be a Polish space, and let $p\ge 1$. For all $\mu,\nu\in\mathcal P_p(\mathcal E)$
\[
\mathcal W_p(\mu,\nu)^p=\sup_{a(x)+b(y)\le d(x,y)^p}\left(\int_\mathcal E a(x)\mu(dx)+\int_\mathcal E b(y)\nu(dy)\right)\,.
\]
In this formula, it is equivalent to assume that the functions $a$ and $b$ belong to $C_b(\mathcal E)$, or that $a\in L^1(\mathcal E,\mu)$ while $b\in L^1(\mathcal E,\nu)$, in which case $a(x)+b(y)\le d(x,y)^p$ 
holds $\mu\otimes\nu$-a.e.. (See for instance Theorem 1.3 in \cite{VillaniAMS}.)

Henceforth, in the latter case, we shall systematically normalize the pair $(a,b)$ as follows: pick $a,b:\,\mathcal E\to\mathbf R\cup\{-\infty\}$ to be measurable representatives of the corresponding elements in 
$L^1(\mathcal E,\mu)$ and $L^1(\mathcal E,\nu)$ resp., together with a $\mu$-negligible set $M\subset\mathcal E$ and a $\nu$-negligible set $N\subset\mathcal E$ such that $a(x)+b(y)\le d(x,y)^p$ holds
for all $x\in M^c$ and all $y\in N^c$. Modifying $a$ and $b$ so that $a(x)=-\infty$ for all $x\in M$ and $b(y)=-\infty$ for all $y\in N$, we obtain in this way new measurable representatives of the same elements
$L^1(\mathcal E,\mu)$ and $L^1(\mathcal E,\nu)$ as before, such that the inequality $a(x)+b(y)\le d(x,y)^p$ holds for all $x,y\in\mathcal E$.

\smallskip
It will be convenient to define a notion of $p$-Legendre transform, as follows: for $f:\,\mathcal E\to\mathbf R\cup\{-\infty\}$, not identically equal to $-\infty$, set
\[
f^{[p*]}(y):=\inf_{x\in\mathcal E}\left(d(x,y)^p-f(x)\right)\,,\quad y\in\mathcal E\,.
\]
A function $g:\,\mathcal E\to\mathbf R\cup\{-\infty\}$ is said to be $d^p$-concave, if it is of the form $g=f^{[p*]}$ for some $f:\,\mathcal E\to\mathbf R\cup\{-\infty\}$ that is not identically equal to $-\infty$.

\smallskip
\noindent 
\textbf{Optimal Kantorovich Potentials.} Under the same assumptions as in the Kantorovich duality theorem for $\mathcal W_p$, there exists a pair $(a,b)$ of $p$-Legendre conjugate functions --- i.e. $b=a^{[p*]}$
and $a=b^{[p*]}$ --- such that $a\in L^1(\mathcal E,\mu)$ and $b\in L^1(\mathcal E,\nu)$, and 
\[
\mathcal W_p(\mu,\nu)^p=\int_\mathcal E a(x)\mu(dx)+\int_\mathcal E b(y)\nu(dy)\,.
\]

(See Remark 1.12, the double convexification trick (2.10), Remark 2.2 and Exercise 2.36 in \cite{VillaniAMS}.)

\smallskip
There are two particular cases of special interest, corresponding to $p=1$ and $p=2$. 

\smallskip
\noindent
\textbf{Kantorovich-Rubinstein Duality Theorem.} Let $(\mathcal E,d)$ be a Polish space with metric $d$. Then, for all $\mu,\nu\in\mathcal P_1(\mathcal E)$, it holds
\[
\mathcal W_1(\mu,\nu)=\sup_{\text{Lip}(\phi)\le 1}\left|\int_\mathcal E\phi(z)\mu(dz)-\int_\mathcal E\phi(z)\nu(dz)\right|\,.
\]

\smallskip
(Indeed, one can check that $(\phi^{[1*]})^{[1*]}=-\phi^{[1*]}$ is a contraction on $\mathcal E$: see Theorem 1.14 and its proof in \cite{VillaniAMS}.)

\smallskip
In the case $p=2$, assuming that $\mathcal E$ is a Euclidean space and $d(x,y)=|x-y|$ is its Euclidean metric, the notion of $2$-Legendre duality is easily reduced to the classical notion of Legendre transform (as
in \S 26 of \cite{Rockaf}, for instance). Indeed
\[
g(y)=\inf_{x\in\mathcal E}(|x-y|^2-f(x))=|y|^2+\inf_{x\in\mathcal E}(|x|^2-2x\cdot y-f(x))
\]
if and only if
\[
\tfrac12(|y|^2-g(y))=\sup_{x\in\mathcal E}(x\cdot y-\tfrac12(|x|^2-f(x)))\,.
\]
In other words, defining $F(x):=\tfrac12(|x|^2-f(x))$ and $G(y):=\tfrac12(|y|^2-g(y))$, it holds
\[
g=f^{[2*]}\iff G=F^*\text{ (the Legendre transform of $F$). }
\]


\section{A Duality-Based Proof of the Triangle Inequality for $\mathcal W_2$}


In this section, we explain how to use the Kantorovich duality to prove that
\[
\mathcal W_2(\lambda,\nu)\le\mathcal W_2(\lambda,\mu)+\mathcal W_2(\mu,\nu)
\]
for all $\lambda,\mu,\nu\in\mathcal P_2(\mathcal E)$. This first result is a warm-up for the duality-based proof of the triangle inequality for $\mathcal W_p$ in the case of an arbitrary exponent $p\in[1,+\infty)$. We recall that 
$(\mathcal E,d)$ is a Polish space with metric denoted by $d$. 

Pick $2$-Legendre conjugate, optimal Kantorovich potentials for $\mathcal W_2(\lambda,\nu)$, denoted by $\alpha$ and $\gamma$. Thus
$$
\alpha(x)=\inf_{z\in\mathcal E}(d(x,z)^2-\gamma(z))\,,\qquad\gamma(z)=\inf_{x\in\mathcal E}(d(x,z)^2-\alpha(x))\,;
$$
besides $\alpha\in L^1(\mathcal E,\lambda)$ while $\gamma\in L^1(\mathcal E,\nu)$, and
$$
\mathcal W_2(\lambda,\nu)^2=\int_{\mathcal E}\alpha(x)\lambda(dx)+\int_{\mathcal E}\gamma(z)\nu(dz)\,.
$$

For each $\eta>0$, set
\begin{equation}\label{DefBeta}
\beta_\eta(y):=(1+\tfrac1\eta)\inf_{z\in\mathcal E}\left(d(y,z)^2-\tfrac{\gamma(z)}{1+\frac1\eta}\right)\,.
\end{equation}
Obviously, the function 
\[
y\mapsto\beta_\eta(y)/(1+\tfrac1\eta)
\]
is $d^2$-concave, since $\gamma\in L^1(\mathcal E,\nu)$ is $\nu$-a.e. finite.

\noindent
\textbf{Lemma 1.} Under the assumptions above, it holds
$$
\alpha(x)-\beta_\eta(y)\le(1+\eta)d(x,y)^2\,,\qquad\text{ for all }x,y\in\mathcal E\,.
$$

\smallskip
Taking Lemma 1 for granted, we conclude the proof of the triangle inequality for $\mathcal W_2$.

\begin{proof}[Proof of the triangle inequality for $\mathcal W_2$] First observe that $\beta_\eta\in L^1(\mathcal E;\mu)$. Indeed, 
\[
\begin{aligned}
|\beta_\eta(y)|\le&(1+\tfrac1\eta)d(y,z)^2+|\gamma(z)|
\\
\le& (1+\tfrac1\eta)(d(x_0,y)+d(x_0,z))^2+|\gamma(z)|
\\
\le&2(1+\tfrac1\eta)(d(x_0,y)^2+d(x_0,z)^2)+|\gamma(z)|\,,
\end{aligned}
\]
where $x_0$ is any point chosen in $\mathcal E$. Integrating both sides of this inequality with the measure $\mu\otimes\nu$, we find that
\[
\int_\mathcal E|\beta_\eta(y)|\mu(dy)\le 2(1+\tfrac1\eta)\left(\int_\mathcal Ed(x_0,y)^2\mu(dy)+\int_\mathcal Ed(x_0,z)^2\nu(dz)\right)+\int_\mathcal E|\gamma(z)|\nu(dz)
\]
so that $\beta_\eta\in L^1(\mathcal E,\mu)$ since $\mu,\nu\in\mathcal P_2(\mathcal E)$ and $\gamma\in L^1(\mathcal E,\nu)$.

Since $\beta_\eta\in L^1(\mathcal E,\mu)$ and $\gamma\in L^1(\mathcal E,\nu)$ satisfy
\[
\tfrac{\beta_\eta(y)}{1+\frac1\eta}+\tfrac{\gamma(z)}{1+\frac1\eta}\le d(y,z)^2\,,\qquad y,z\in\mathcal E
\]
according to the definition \eqref{DefBeta} of $\beta_\eta$, the Kantorovich Duality Theorem implies that
\begin{equation}\label{Wmunu}
\tfrac1{1+\frac1\eta}\int_\mathcal E\beta_\eta(y)\mu(dy)+\tfrac1{1+\frac1\eta}\int_\mathcal E\gamma(z)\nu(dz)\le\mathcal W_2(\mu,\nu)^2\,.
\end{equation}

On the other hand, the fact that $\alpha\in L^1(\mathcal E,\lambda)$ and $\beta_\eta\in L^1(\mathcal E,\mu)$, together with the inequality in Lemma 1 and the Kantorovich Duality Theorem, imply that
\begin{equation}\label{Wlamu}
\tfrac1{1+\eta}\int_\mathcal E\alpha(x)\lambda(dx)-\tfrac1{1+\eta}\int_\mathcal E\beta_\eta(y)\mu(dy)\le\mathcal W_2(\lambda,\mu)^2\,.
\end{equation}

Multiplying both sides of \eqref{Wmunu} by $1+\tfrac1\eta$ and both sides of \eqref{Wlamu} by $1+\eta$, and adding each side of the resulting inequalities, we see that
\[
\begin{aligned}
\mathcal W_2(\lambda,\nu)^2=\int_\mathcal E\alpha(x)\lambda(dx)+\int_\mathcal E\gamma(z)\nu(dz)
\\
=\int_\mathcal E\alpha(x)\lambda(dx)-\int_\mathcal E\beta_\eta(y)\mu(dy)+\int_\mathcal E\beta_\eta(y)\mu(dy)+\int_\mathcal E\gamma(z)\nu(dz)
\\
\le(1+\eta)\mathcal W_2(\lambda,\mu)^2+(1+\tfrac1\eta)\mathcal W_2(\mu,\nu)^2&\,.
\end{aligned}
\]
Assuming that $\lambda\not=\mu$, so that $\mathcal W_2(\lambda,\mu)>0$, we pick $\eta:=\mathcal W_2(\mu,\nu)/\mathcal W_2(\lambda,\mu)$ to find that
\[
\begin{aligned}
\mathcal W_2(\lambda,\nu)^2\le&\mathcal W_2(\lambda,\mu)^2+\mathcal W_2(\mu,\nu)^2+2\mathcal W_2(\lambda,\mu)\mathcal W_2(\mu,\nu)
\\
=&(\mathcal W_2(\lambda,\mu)+\mathcal W_2(\mu,\nu))^2\,.
\end{aligned}
\]
Otherwise, $\lambda=\mu$ and the triangle inequality is trivial.
\end{proof}

It only remains to prove Lemma 1.

\begin{proof}[Proof of Lemma 1.] We seek to bound
\[
\begin{aligned}
\alpha(x)-\beta_\eta(y)=\inf_{z\in\mathcal E}(d(x,z)^2-\gamma(z))-(1+\tfrac1\eta)\inf_{z\in\mathcal E}\left(d(y,z)^2-\tfrac{\gamma(z)}{1+\frac1\eta}\right)
\\
=\inf_{z\in\mathcal E}(d(x,z)^2-\gamma(z))-\inf_{z\in\mathcal E}\left((1+\tfrac1\eta)d(y,z)^2-\gamma(z)\right)
\end{aligned}
\]
Let $\eps>0$; there exists $z_\eps\in\mathbf R^d$ such that
\[
\begin{aligned}
\inf_{z\in\mathcal E}\left((1+\tfrac1\eta)d(y,z)^2-\gamma(z)\right)\le&\left((1+\tfrac1\eta)d(y,z_\eps)^2-\gamma(z_\eps)\right)
\\
<&\inf_{z\in\mathcal E}\left((1+\tfrac1\eta)d(y,z)^2-\gamma(z)\right)+\eps\,.
\end{aligned}
\]
Then
\[
\begin{aligned}
\alpha(x)-\beta_\eta(y)=&\inf_{z\in\mathcal E}(d(x,z)^2-\gamma(z))-\inf_{z\in\mathcal E}\left((1+\tfrac1\eta)d(y,z)^2-\gamma(z)\right)
\\
\le&(d(x,z_\eps)^2-\gamma(z_\eps))-\inf_{z\in\mathcal E}\left((1+\tfrac1\eta)d(y,z)^2-\gamma(z)\right)
\\
<&(d(x,z_\eps)^2-\gamma(z_\eps))-\left((1+\tfrac1\eta)d(y,z_\eps)^2-\gamma(z_\eps)\right)+\eps
\\
=&d(x,z_\eps)^2-(1+\tfrac1\eta)d(y,z_\eps)^2+\eps\,.
\end{aligned}
\]
But, for each $\eta>0$, it holds
$$
\begin{aligned}
d(x,z_\eps)^2\le&(d(x,y)+d(y,z_\eps))^2
\\
=&d(x,y)^2+d(y,z_\eps)^2+2d(x,y)d(y,z_\eps)
\\
\le&d(x,y)^2+d(y,z_\eps)^2+\eta d(x,y)^2+\tfrac1\eta d(y,z_\eps)^2
\\
=&(1+\eta)d(x,y)^2+(1+\tfrac1\eta)d(y,z_\eps)^2\,.
\end{aligned}
$$
With the preceding inequality, we conclude that
\[
\alpha(x)-\beta_\eta(y)\le(1+\eta)d(x,y)^2+\eps\,,
\]
and the desired inequality follows from letting $\eps\to 0^+$.
\end{proof}


\section{A Duality-Based Proof of the Triangle Inequality\\ for $\mathcal W_p$ with $1\le p<\infty$ and $p\not=2$}


In this section, we use the Kantorovich duality to prove that
\begin{equation}\label{Triang-p}
\mathcal W_p(\lambda,\nu)\le\mathcal W_p(\lambda,\mu)+\mathcal W_p(\mu,\nu)
\end{equation}
for all $p\ge 1$ and all $\lambda,\mu,\nu\in\mathcal P_2(\mathcal E)$. We recall that $(\mathcal E,d)$ is a Polish space with metric denoted by $d$.

The case $p=1$ follows immediately from the formula for $\mathcal W_1$ in the Kantorovich-Rubinstein Duality Theorem. Henceforth, we assume therefore that
\[
p>1\quad\text{ and }\quad p\not=2\,.
\]

A careful inspection of the duality-based proof of the triangle inequality for $\mathcal W_2$ shows the importance of the inequality
\begin{equation}
\label{Ineq-2}
(X+Y)^2\le(1+\eta)X^2+(1+\tfrac1\eta)Y^2
\end{equation}
for all $X,Y\ge 0$ and all $\eta>0$. 

Our first task is therefore to seek a function $(0,+\infty)\ni\eta\mapsto f(\eta)\in(0,+\infty)$ such that
\[
(X+Y)^p\le(1+\eta)X^p+(1+f(\eta))Y^p\,,\qquad X,Y\ge 0\,,\quad\eta>0\,.
\]
Obviously, only the case $X,Y>0$ is of interest, so that, by homogeneity, this boils down to finding $f$ such that
\[
(Z^{1/p}+1)^p\le(1+\eta)Z+1+f(\eta)\,,\qquad Z,\eta>0\,,
\]
where $Z=X^p/Y^p$.

Equivalently, the optimal $f(\eta)$ is found to be given by the formula
\[
f(\eta):=\sup_{Z>0}(-\eta Z-1-Z+(1+Z^{1/p})^p)\,.
\]
One easily checks that the function $(0,+\infty)\ni Z\mapsto(1+Z^{1/p})^p\in(0,+\infty)$ is concave for $p>1$, since
\[
\tfrac{d}{dZ}(1+Z^{1/p})^p=p(1+Z^{1/p})^{p-1}\tfrac1p Z^{\frac1p-1}=(Z^{-1/p}+1)^{p-1}
\]
defines a decreasing bijection from $(0,+\infty)$ to itself. Hence there exists a unique critical value of $Z>0$ such that
\[
\tfrac{d}{dZ}(-\eta Z-1-Z+(1+Z^{1/p})^p)=-(\eta+1)+(Z^{-1/p}+1)^{p-1}=0\,,
\]
which is
\[
Z^{1/p}:=\frac1{(\eta+1)^{1/(p-1)}-1}\,,
\]
and $f$ is given by the formula
\[
\begin{aligned}
f(\eta):=&\left(1+\frac1{(\eta+1)^{1/(p-1)}-1}\right)^p-1-\frac{\eta+1}{((\eta+1)^{1/(p-1)}-1)^p}
\\
=&\left(\frac{(\eta+1)^{1/(p-1)}}{(\eta+1)^{1/(p-1)}-1}\right)^p-1-\frac{\eta+1}{((\eta+1)^{1/(p-1)}-1)^p}
\\
=&\frac{(\eta+1)^{p/(p-1)}-(\eta+1)}{((\eta+1)^{1/(p-1)}-1)^p}-1\,.
\end{aligned}
\]
Summarizing, we have proved the following lemma.

\medskip
\noindent
\textbf{Lemma 2.} For all $X,Y\ge 0$ and all $\eta>0$, it holds
\[
(X+Y)^p\le(1+\eta)X^p+\frac{(\eta+1)^{p/(p-1)}-(\eta+1)}{((\eta+1)^{1/(p-1)}-1)^p}Y^p\,.
\]

\smallskip
(One easily checks that, in the case $p=2$, 
\[
f(\eta)=\frac{(\eta+1)^2-(\eta+1)}{((\eta+1)-1)^2}-1=\frac{\eta^2+\eta}{\eta^2}-1=\frac1\eta
\]
so that the inequality in Lemma 2 coincides with \eqref{Ineq-2}.)

\smallskip
Next we prove \eqref{Triang-p}. Pick $p$-Legendre conjugate, optimal Kantorovich potentials for $\mathcal W_p(\lambda,\nu)$, denoted by $\alpha$ and $\gamma$ as in the preceding section. Thus
$$
\alpha(x)=\inf_{z\in\mathcal E}(d(x,z)^p-\gamma(z))\,,\qquad\gamma(z)=\inf_{x\in\mathcal E}(d(x,z)^p-\alpha(x))\,;
$$
besides $\alpha\in L^1(\mathcal E,\lambda)$ while $\gamma\in L^1(\mathcal E,\nu)$, and
$$
\mathcal W_p(\lambda,\nu)^p=\int_{\mathcal E}\alpha(x)\lambda(dx)+\int_{\mathcal E}\gamma(z)\nu(dz)\,.
$$

For each $\eta>0$, set
\begin{equation}\label{DefBeta-p}
\beta_\eta(y):=\frac{(\eta+1)^{p/(p-1)}-(\eta+1)}{((\eta+1)^{1/(p-1)}-1)^p}\inf_{z\in\mathcal E}\left(d(y,z)^p-\frac{((\eta+1)^{1/(p-1)}-1)^p}{(\eta+1)^{p/(p-1)}-(\eta+1)}\gamma(z)\right)\,.
\end{equation}
Obviously, the function 
\[
y\mapsto\frac{((\eta+1)^{1/(p-1)}-1)^p}{(\eta+1)^{p/(p-1)}-(\eta+1)}\beta_\eta(y)
\]
is $d^p$-concave, since $\gamma\in L^1(\mathcal E,\nu)$ is $\nu$-a.e. finite.

\medskip
\noindent
\textbf{Lemma 3.} Under the assumptions above, it holds
\[
\alpha(x)-\beta_\eta(y)\le(1+\eta)d(x,y)^p\,,\qquad x,y\in\mathcal E\,.
\]

\begin{proof}[Proof of Lemma 3] For each $\eps>0$, there exists $z_\eps\in\mathcal E$ such that
\[
\begin{aligned}
\inf_{z\in\mathcal E}&\left(\frac{(\eta+1)^{p/(p-1)}-(\eta+1)}{((\eta+1)^{1/(p-1)}-1)^p}d(y,z)^p-\gamma(z)\right)
\\
&\le\left(\frac{(\eta+1)^{p/(p-1)}-(\eta+1)}{((\eta+1)^{1/(p-1)}-1)^p}d(y,z_\eps)^p-\gamma(z_\eps)\right)
\\
&<\inf_{z\in\mathcal E}\left(\frac{(\eta+1)^{p/(p-1)}-(\eta+1)}{((\eta+1)^{1/(p-1)}-1)^p}d(y,z)^p-\gamma(z)\right)+\eps
\end{aligned}
\]
Thus
\[
\begin{aligned}
\alpha(x)-\beta_\eta(y)=&\inf_{z\in\mathcal E}(d(x,z)^p-\gamma(z))-\inf_{z\in\mathcal E}\left(\tfrac{(\eta+1)^{p/(p-1)}-(\eta+1)}{((\eta+1)^{1/(p-1)}-1)^p}d(y,z)^p-\gamma(z)\right)
\\
<&(d(x,z_\eps)^p-\gamma(z_\eps))-\left(\tfrac{(\eta+1)^{p/(p-1)}-(\eta+1)}{((\eta+1)^{1/(p-1)}-1)^p}d(y,z_\eps)^p-\gamma(z_\eps)\right)+\eps
\\
\le&(d(x,y)+d(y,z_\eps))^p-\tfrac{(\eta+1)^{p/(p-1)}-(\eta+1)}{((\eta+1)^{1/(p-1)}-1)^p}d(y,z_\eps)^p+\eps
\\
\le&(1+\eta)d(x,y)^p+\eps
\end{aligned}
\]
by Lemma 2. We conclude by letting $\eps\to 0^+$.
\end{proof}

\smallskip
Assume that $\mathcal W_p(\lambda,\mu)\mathcal W_p(\mu,\nu)\not=0$. First observe that
\[
|\beta_\eta(y)|\le\tfrac{(\eta+1)^{p/(p-1)}-(\eta+1)}{((\eta+1)^{1/(p-1)}-1)^p}(d(x_0,y)+d(x_0,z))^p+|\gamma(z)|\,,
\]
and integrating both sides of this inequality with the measure $\mu\otimes\nu$ shows that
\[
\begin{aligned}
\int_\mathcal E|\beta_\eta(y)|\mu(dy)\le&2^{p-1}\tfrac{(\eta+1)^{p/(p-1)}-(\eta+1)}{((\eta+1)^{1/(p-1)}-1)^p}\left(\int_\mathcal Ed(x_0,y)^p\mu(dy)+\int_\mathcal Ed(x_0,z)^p\nu(dy)\right)
\\
&+\int_\mathcal E|\gamma(z)|\nu(dz)<+\infty\,,
\end{aligned}
\]
since $\mu,\nu\in\mathcal P_p(\mathcal E)$ and $\gamma\in L^1(\mathcal E,\nu)$. Hence $\beta_\eta\in L^1(\mathcal E,\mu)$ and the inequality in Lemma 3 implies that
\[
\int_\mathcal E\alpha(x)\lambda(dx)-\int_\mathcal E\beta_\eta(y)\mu(dy)\le(1+\eta)\mathcal W_p(\lambda,\mu)^p
\]
by Kantorovich duality. On the other hand, the definition \eqref{DefBeta-p} implies that
\[
\beta_\eta(y)+\gamma(z)\le\tfrac{(\eta+1)^{p/(p-1)}-(\eta+1)}{((\eta+1)^{1/(p-1)}-1)^p}d(y,z)^p\,,\qquad y,z\in\mathcal E\,,
\]
so that
\[
\int_\mathcal E\beta_\eta(y)\mu(dy)+\int_\mathcal E\gamma(z)\nu(dz)\le\tfrac{(\eta+1)^{p/(p-1)}-(\eta+1)}{((\eta+1)^{1/(p-1)}-1)^p}\mathcal W_p(\mu,\nu)^p
\]
again by Kantorovich duality. Hence
\begin{equation}\label{Ineq-eta-p}
\begin{aligned}
\mathcal W_p(\lambda,\nu)^p=&\int_\mathcal E\alpha(x)\lambda(dx)+\int_\mathcal E\gamma(z)\nu(dz)
\\
\le&(1+\eta)\mathcal W_p(\lambda,\mu)^p+\tfrac{(\eta+1)^{p/(p-1)}-(\eta+1)}{((\eta+1)^{1/(p-1)}-1)^p}\mathcal W_p(\mu,\nu)^p\,,
\end{aligned}
\end{equation}
and this inequality holds for each $\eta>0$. Choose
\[
\eta+1:=(Z^{-1/p}+1)^{p-1}\,,\quad\text{ with }\quad Z:=\mathcal W_p(\lambda,\mu)^p/\mathcal W_p(\mu,\nu)^p\,.
\]
Then
\[
\begin{aligned}
(1+\eta)Z+\tfrac{(\eta+1)^{p/(p-1)}-(\eta+1)}{((\eta+1)^{1/(p-1)}-1)^p}=&(Z^{-1/p}+1)^{p-1}Z+\tfrac{(Z^{-1/p}+1)^p-(Z^{-1/p}+1)^{p-1}}{(Z^{-1/p}+1-1)^p}
\\
=&(1+Z^{1/p})^{p-1}Z^{1/p}+\tfrac{(Z^{-1/p}+1-1)(Z^{-1/p}+1)^{p-1}}{(Z^{-1/p}+1-1)^p}
\\
=&(1+Z^{1/p})^{p-1}Z^{1/p}+\tfrac{Z^{-1/p}(Z^{-1/p}+1)^{p-1}}{Z^{-1}}
\\
=&(1+Z^{1/p})^{p-1}Z^{1/p}+Z^{(p-1)/p}(Z^{-1/p}+1)^{p-1}
\\
=&(1+Z^{1/p})^{p-1}Z^{1/p}+(1+Z^{1/p})^{p-1}
\\
=&(1+Z^{1/p})^p\,.
\end{aligned}
\]
In other words, with this choice of $\eta$ and $Z$, one finds that
\[
(1+\eta)\frac{\mathcal W_p(\lambda,\mu)^p}{\mathcal W_p(\mu,\nu)^p}+\frac{(\eta+1)^{p/(p-1)}-(\eta+1)}{((\eta+1)^{1/(p-1)}-1)^p}=\left(1+\frac{\mathcal W_p(\lambda,\mu)}{\mathcal W_p(\mu,\nu)}\right)^p\,.
\]
Multiplying both sides of this identity by $\mathcal W_p(\mu,\nu)^p$, one arrives at the identity
\[
(1+\eta)\mathcal W_p(\lambda,\mu)^p+\frac{(\eta+1)^{p/(p-1)}-(\eta+1)}{((\eta+1)^{1/(p-1)}-1)^p}\mathcal W_p(\mu,\nu)^p=(\mathcal W_p(\lambda,\mu)+\mathcal W_p(\mu,\nu))^p
\]
where $\eta$ is chosen as follows:
\[
\eta:=\left(\tfrac{\mathcal W_p(\mu,\nu)}{\mathcal W_p(\lambda,\mu)}+1\right)^{p-1}-1\,.
\] 
Inserting this value of $\eta$ in the right-hand side of \eqref{Ineq-eta-p} and using the identity above shows that
\[
\mathcal W_p(\lambda,\nu)^p\le(\mathcal W_p(\lambda,\mu)+\mathcal W_p(\mu,\nu))^p
\]
from which the triangle inequality immediately follows in the special  case where 
\[
\mathcal W_p(\lambda,\mu)\mathcal W_p(\mu,\nu)\not=0\,.
\]

Otherwise, if one of the distances $\mathcal W_p(\lambda,\mu)$ or $\mathcal W_p(\mu,\nu)$ is equal to $0$, the triangle inequality is trivial.


\section{Conclusion}


The proof of the triangle inequality for the Wasserstein distances presented in this short note is (perhaps) the most ``elementary'' proof --- which does not mean that it is the shortest --- in the sense that it does 
not use any information about optimal couplings. (For instance, it does not use the existence of an optimal transport map, as in the Brenier theorem, i.e. Theorem 2.12 (ii) in \cite{VillaniAMS}.) Also, it does not 
involve any nontrivial manipulation on optimal couplings (as in the ``glueing'' procedure described in Lemma 7.6 of \cite{VillaniAMS}). The only ingredients in this proof are (a) Kantorovich duality, and (b) the 
inequality of Lemma 2, which boils down to computing the Legendre transform of a real-valued function on the half-line. The proof of Lemma 3 is based on (b)  --- just as Lemma 1 is based on the elementary 
inequality $2ab\le\eta a^2+\tfrac1\eta b^2$ for all $a,b,\eta>0$ --- and on the characterization of the infimum of a function as its larger lower bound. Since there is a Kantorovich-type duality for the quantum
analogue of $\mathcal W_2$ defined in \cite{FGMouPaul} (see \cite{CaglioFGPaul}), it seems that the validity of the triangle inequality for this quantum Wasserstein (pseudo)metric boils down to the existence
of a quantum analogue of Lemma 1. Since D. Serre's example mentioned in the introduction rules out the possibility of proving the triangle inequality for the quantum analogue of $\mathcal W_2$ defined in
\cite{FGMouPaul} by glueing quantum couplings as explained on p. 28 of \cite{FGPaulJFA}, we hope that the approach to the triangle inequality presented here can shed light on the quantum case.



\end{document}